\documentclass{lms}
\usepackage{epsfig}
\usepackage{amssymb}
\begin{document}

\date{}

\title{ERRATA TO `THE LIMITING CURVE OF JARN\'{I}K'S POLYGONS'}
\author{Jovi\v{s}a \v{Z}uni\'c}

\classno{11P21 (primary); 52C05 (secondary)}

\maketitle

\begin{abstract} 
\begin{itemize} 
\item[] In this note we point out that the limit shape of Jarn\'{i}k's polygonal curve \cite{jarnik} is the circle, not a curve 
consisting of arcs of parabolas as it has been stated in the main result of \cite{paper-AMS} (Section 2). The correct 
result has been established and proven in \cite{paper-AA}. 
%\\[0.2cm] Both papers \cite{paper-AMS} and \cite{paper-AA} are use the idea and technique from  \cite{paper-DM}.
\end{itemize}
\end{abstract}
%\vspace*{0.3cm}
\section{Introduction}
In order to establish the maximum number $f(s)$ of lattice points on a strictly convex curve 
of length at most $s$, Jarn\'{i}k \cite{jarnik} has considered a family of convex lattice polygons.
Vertices of these polygons lie on strictly convex curves which optimize $f(s)$, in an asymptotic sense.
It has been shown \cite{paper-AA} that there is a limit shape of such polygons. More precisely, as the 
number of vertices of these polygons increase and after a proper scaling (e.g., for a factor equal to the 
perimeter of the related polygon) these polygons converge to a circle.  
The result has been derived starting from the family of the convex lattice polygons whose edge 
slopes $q/a$ belong to the set $V_t$, defined as follows (for an arbitrary $t>1$):
\begin{equation}
V_t = \left\{(q, a) \in \mathbf{Z}^2 \ : \  \gcd(q, a) = 1; \ \ \sqrt{q^2 + a^2} \leq t \right\}.
\end{equation}

Another family of convex lattice polygons has been considered in \cite{paper-AMS}. The edge slopes 
$q/a$ of lattice polygons from this family make the set $V_Q,$ defined for all $Q>1$, as follows
\begin{equation}
V_Q = \left\{(q, a) \in \mathbf{Z}^2 \ : \  \gcd(q, a) = 1; \ \ \max\{|a|, |q|\} \leq Q \right\}.
\end{equation}
The polygons $V_Q$ do not relate to the Jarn\'{i}k's polygons (i.e., do not optimize the number of the polygon
vertices with respect to the Euclidean perimeter of the polygon), as it has been stated in \cite{paper-AMS}. 
Indeed, the observation of $V_Q$, as given in  \cite{paper-AMS}, does not lead to a conclusion about the 
number of lattice points on a strictly convex curve with respect to the Euclidean length of this curve.

Thus, the statement, from \cite{paper-AMS},
that the limiting curve of Jarn\'{i}k's polygons consist of parabolic arcs is not true.
The limit shape of such polygons is a circle, as it has been shown in \cite{paper-AA}. The result 
from \cite{paper-AMS} can be considered in sense of the maximal number of lattice points, on a strictly convex
curve, with respect to the curve perimeter taken in sense of $l_\infty$ distance, rather than in sense of the $l_2$ (i.e. Euclidean)
distance.

It is also worth mentioning the following: The some of the results  from \cite{paper-AMS} (see Section 2) 
are derived following the idea and technique from \cite{paper-AA} and \cite{paper-DM}, but this was not referenced properly;
The term `limit shape' has been used by the others (see \cite{vershik}, for an example) but in different sense than 
it has been done in \cite{paper-AMS}, \cite{paper-AA}, and \cite{paper-DM}.

\affiliationone{
Jovi\v{s}a \v{Z}uni\'{c}  \\ 
Mathematical Institute \\ Serbian Academy of Sciences and Arts\\ 
Belgrade, Serbia
\email{jovisa\_zunic@mi.sanu.ac.rs}}

\end{document}